\documentclass[12pt,leqno]{article}

\newtheorem{theorem}{Theorem}[section]
\newtheorem{lemma}[theorem]{Lemma}
\newtheorem{proposition}[theorem]{Proposition}
\newtheorem{corollary}[theorem]{Corollary}
\newtheorem{definition}[theorem]{Definition}
\newcommand{\pf}{\noindent{\bf Proof. }}
\title{\Large \bf Compactness and an approximation property related to an operator ideal}
\author{Anil Kumar Karn$^{\dagger}$ and Deba Prasad Sinha$^{\ddagger}$}

\date{}
\begin{document}

\maketitle

\thispagestyle{empty}

\begin{abstract}
For an operator ideal $\mathcal A$, we study the composition operator ideals ${\mathcal A}\circ{\mathcal K}$,
${\mathcal K}\circ{\mathcal A}$ and ${\mathcal K}\circ{\mathcal A}\circ{\mathcal K}$, where $\mathcal K$ is the
ideal of compact operators. We introduce a notion of an $\mathcal A$-approximation property on a Banach space
and characterise it in terms of the density of finite rank operators in ${\mathcal A}\circ{\mathcal K}$ and
${\mathcal K}\circ{\mathcal A}$.

We propose the notions of $\ell _{\infty}$-extension and $\ell _1$-lifting properties for an operator ideal
$\mathcal A$ and study ${\mathcal A}\circ{\mathcal K}$, ${\mathcal K}\circ{\mathcal A}$ and
the $\mathcal A$-approximation property where $\mathcal A$ is injective or surjective and/or
with the $\ell _{\infty}$-extension or $\ell _1$-lifting property.  In particular, we show that if $\mathcal A$ is an injective
operator ideal with the $\ell _\infty$-extension property, then we have:
{\item{(a)} $X$ has the $\mathcal A$-approximation property if and only if
$({\mathcal A}^{min})^{inj}(Y,X)={\mathcal A}^{min}(Y,X)$,
for all Banach spaces $Y$.
\item{(b)} The dual space $X^*$ has the $\mathcal A$-approximation property if and only if
$(({\mathcal A}^{dual})^{min})^{sur}(X,Y)=({\mathcal A}^{dual})^{min}(X,Y)$,
for all Banach spaces $Y$.}For an operator ideal $\mathcal A$, we study the composition operator ideals ${\mathcal A}\circ{\mathcal K}$,
\end{abstract}

\begin{description}
\item[$\dagger$] School of Mathematical Sciences, National Institute of Science Education and Research, Institute of Physics Campus, P.O. Sainik School, Bhubaneswar-751005, India; e-mail: {\texttt anilkarn@niser.ac.in}
\item[$\ddagger$] (Lately at) Department of Mathematics, Dyal Singh College (University of Delhi), Lodi Road, New Delhi 110003, India; e-mail: {\texttt sinha.deba@gmail.com}
\item [2010  Mathematics Subject Classification:] Primary 46B50; Secondary 46B20, 46B28, 47B07
\item [Keywords and phrases:] approximation property, kernel procedure, injective operator ideal, surjective operator ideal, $\ell _\infty$-extension property, $\ell _1$-lifting property, $\mathcal A$-approximation property.
\end{description}

\section{Introduction}
It is well known that a Banach space $Y$ has the approximation property if and only if,
$\overline{\mathcal F}(X, Y)
= {\mathcal K}(X, Y)$ for all Banach spaces $X$. Similarly, the dual $X^*$ of a Banach
space $X$ has the approximation property if and only if $\overline{\mathcal F}(X, Y)
= {\mathcal K}(X, Y)$ for all Banach spaces $Y$. However, in general for a pair
of Banach spaces $X$ and $Y$,
$\overline{\mathcal F}(X, Y) \neq {\mathcal K}(X, Y)$, whereas
$\overline{\mathcal F}(X, Y)^{inj}
= {\mathcal K}(X, Y) = \overline{\mathcal F}(X, Y)^{sur}$.
In the language of operator ideals $\overline{\mathcal F}^{inj}
= {\mathcal K} = \overline{\mathcal F}^{sur}$.   In this language it may be stated that a
Banach space $X$ has the approximation property if and only if $$\overline{\mathcal F}^{inj}(Y,X)
= \overline{\mathcal F}(Y,X)$$ for all Banach spaces $Y$ and the dual space $X^*$ has the
approximation property if and only if $$\overline{\mathcal F}^{sur}(X,Y) = \overline{\mathcal F}(X,Y)$$
for all Banach spaces $Y$.

In the papers \cite{SK1, SK2} , the authors introduced a class of operator ideals ${\mathcal K}_p$, $(1 \leq p < \infty)$ of compact
operators whose adjoints factor through specific subspaces of $l_p$ and showed that
$({\Pi _p}^{min})^{inj} = {\mathcal K}_p^{dual} = {\Pi _p} \circ {\mathcal K}$ and
$(({\Pi _p}^{dual})^{min})^{sur} = {\mathcal K}_p^{dual \hspace{1mm}dual} =
{\mathcal K} \circ {\Pi _p}^{dual}$,
where $\Pi _p$ is the operator ideal of $p$-summing operators. (In the limiting case
${\mathcal B}\circ{\mathcal K} = {\mathcal K}\circ{\mathcal B} = {\mathcal K}$ holds
trivially and also we have ${\mathcal B}^{min} = \overline{\mathcal F}$.) Further, in \cite{SK2} they introduced a
notion of the approximation property of {\it type} $p$ (related to the operator ideal $\Pi _p$) and proved that
a Banach space $X$ has the approximation property of {\it type} $p$ if and only if
$$(\Pi _p^{min})^{inj}(Y,X) = \Pi _p^{min}(Y,X)$$ for all Banach spaces $Y$. Similarly, the dual
space $X^*$ has the approximation property of {\it type} $p$ if and only if
$$((\Pi _p^{dual})^{min})^{sur}(X,Y) = (\Pi _p^{dual})^{min}(X,Y)$$ for all Banach spaces $Y$.

The above discussions brings to sharp relief the need for the following:
\indent (I) To study the composition operator
ideals ${\mathcal A}\circ{\mathcal K}$ and ${\mathcal K}\circ{\mathcal A}$, and also,
in general, ${\mathcal K}\circ{\mathcal A}\circ{\mathcal K}$ which we denote by
${\mathcal A}^{com}$ as the compact level objects related to an operator ideal $\mathcal A$, and
\newline\indent (II)To introduce an approximation property related to an operator ideal $\mathcal A$ so
as to extend the above mentioned characterizations of the approximation property to the
operator ideal setting; namely, to study the density of finite rank operators in the relevant compact
level of the operator ideal in the corresponding ideal norm.
\newline\noindent These are the twin objectives of this paper.

In Section 2, we introduce the kernel procedure $com : {\mathcal A} \to {\mathcal K}
\circ{\mathcal A}\circ{\mathcal K}$ for a general operator ideal $\mathcal A$. We
discuss the interplay of this procedure with the standard procedures of operator ideal
theory.

The next two sections are mainly preparatory in nature. Section 3 is of independent
interest wherein we look closely at the
definitions of injective and surjective operator ideals. We show that the composition
of two injective (surjective) operator ideals is again injective (respectively,
surjective).

In Section 4 we introduce two properties, namely, the $l_\infty$-extension
and the $l_1$-lifting properties of operator ideals. These are weaker than the
extension and  lifting
properties that characterise left and right projective operator ideals
\cite[Ex. 20.6]{DF}. We show
that many injective and surjective operator ideals possess these properties. We prove
that if
${\mathcal A}_1$ has the $l_\infty$-extension property, or if ${\mathcal A}_2$
is injective, then
$({\mathcal A}_1
\circ{\mathcal A}_2)^{inj} = {{\mathcal A}_1}^{inj} \circ{{\mathcal A}_2}^{inj}$
and a similar result involving the surjective hull and the $l_1$-lifting property
is also proved. We show that the composition of two
operator ideals both having the $l_\infty$-extension ($l_1$-lifting) property
also has the same property. The duality of these two properties
for an operator ideal $\mathcal A$ and its dual ${\mathcal A}^{dual}$ is also studied.

In Section 5 we investigate the interplay of the kernel procedure {\it `com'}
with the hull procedures {\it `inj'} and {\it `sur'}. We show that
if $\mathcal A$ is left accessible, then
$$ ({\mathcal A}^{com})^{inj} = ({\mathcal A}^{min})^{inj} =
{\mathcal K}\circ{\mathcal A}^{inj}.$$
On the other hand if ${\mathcal A}$ is right accessible and has the
$l_\infty$-extension property, then we have
$$({\mathcal A}^{com})^{inj} = ({\mathcal A}^{min})^{inj} =
{\mathcal A}^{inj}\circ{\mathcal K}.$$
\noindent Dual results for surjective operator ideals and for operator ideals
with the $l_1$-lifting property are also obtained.

In Section 6, we introduce a notion of an approximation property related to an operator ideal
$\mathcal A$, namely the $\mathcal A$-approximation property and prove that a Banach space
$X$ has the $\mathcal A$-approximation property if and only if
$$\overline{{\mathcal F}(X,Y)}^{\alpha _{\kappa}} = {\mathcal A}\circ{\mathcal K}(X,Y),$$ where
$\alpha _{\kappa}$ is the composition ideal norm of ${\mathcal A}\circ{\mathcal K}$.
Similarly, the dual $X^*$ has the $\mathcal A$-approximation property if and only if
$$\overline{{\mathcal F}(X,Y)}^{_{\kappa}{\alpha}^d} = {\mathcal K}\circ{\mathcal A}^{dual}(X,Y),$$
(under certain conditions on $\mathcal A$), where ${_{\kappa}}\alpha ^d$ is the composition
ideal norm of ${\mathcal K}\circ{\mathcal A}^{dual}$.

At the end of the section we study the
$\mathcal A$-approximation property on a Banach space when $\mathcal A$ is injective with the $\ell_\infty$-extension property or
surjective. In fact, we show that if $\mathcal A$ is left accessible injective operator ideal
with the $\ell _{\infty}$-extension property, then a Banach space $X$ has the $\mathcal A$-approximation property
if and only if
$$({\mathcal A}^{min})^{inj}(Y,X) = {\mathcal A}^{min}(Y,X),$$
for all Banach spaces $Y$. The dual space $X^*$ has the $\mathcal A$-approximation property if and only if
$$(({\mathcal A}^{dual})^{min})^{sur}(X,Y) = ({\mathcal A}^{dual})^{min}(X,Y),$$
for all Banach spaces $Y$.

We also show that if  $\mathcal A$ is a right accessible, surjective operator ideal, then $X$ has the
$\mathcal A$-approximation property if and only if
$$({\mathcal A}^{min})^{sur}(Y,X) = {\mathcal A}^{min}(Y,X),$$
for all Banach spaces $Y$. The dual space $X^*$ has the $\mathcal A$-approximation property if and only if
$$(({\mathcal A}^{dual})^{min})^{inj}(X,Y) = ({\mathcal A}^{dual})^{min}(X,Y),$$
for all Banach spaces $Y$.

Before we close the section we take a quick look at some notations. For a Banach space $X$,
$\kappa _X : X \hookrightarrow X^{**}$ is the natural canonical embedding, and in case
$X$ is complemented in its bidual $X^{**}$, $P_X : X^{**} \to X$ shall denote the
resulting projection. For the closed unit ball $B_X$ of $X$, $q_X :l_1(B_X) \to X$
denotes a usual quotient map and $i_X : X \hookrightarrow l_\infty(B_{X^*})$ is the natural
Alouglu embedding. We denote  by $\mathcal B$, $\mathcal K$,
$\overline{\mathcal F}$ and by $\mathit {\Pi}_p$,
${\mathcal I}_p$ and ${\mathcal N}_p$, for $1 \leq p \leq \infty$, the operator ideals
of bounded, compact, approximable, $p$-summing, $p$-integral
and $p$-nuclear operators respectively.

We have avoided routine discussions on norms of operator ideals and that of their compositions
at several places in the body of the paper.

\section{The procedure `com'}

We begin by formally assigning a symbol to ${\mathcal K}\circ{\mathcal A}
\circ{\mathcal K}$ for an operator ideal $\mathcal A$.
\begin{definition}
Let $\mathcal A$ be any quasi normed operator ideal. The composition quasi normed
operator ideal ${\mathcal K}\circ{\mathcal A}\circ{\mathcal K}$ shall be denoted by
${\mathcal A}^{com}$, where $\mathcal K$ is the operator ideal of compact operators.
\end{definition}

It is easy to note that the procedure $com : {\mathcal A} \to {\mathcal A}^{com}$
is monotone. Also as $\mathcal K$ is an idempotent ideal, {\it com} is an idempotent
and a kernel procedure.
In this section we discuss some basic facts regarding the procedure {\it com}
as well as its interplay with other important procedures in the theory of operator ideals .

Let us recall the definition of a minimal kernel ${\mathcal A}^{min}$ of an operator ideal
$\mathcal A$ as the composition operator ideal ${\overline{\mathcal F}}\circ{\mathcal A}
\circ{\overline{\mathcal F}}$,
and the corresponding procedure $min : {\mathcal A} \to {\mathcal A}^{min}$
\cite [Section 4.8]{P1}.
To begin with we consider the interplay of {\it min} and {\it com}. Since,
${\overline{\mathcal F}}\circ{\mathcal K} = {\mathcal K}\circ{\overline{\mathcal F}}
= {\overline{\mathcal F}}$, we have
\begin{proposition}
$({\mathcal A}^{min})^{com} = ({\mathcal A}^{com})^{min} = {\mathcal A}^{min}$.
\end{proposition}

\vskip 8pt plus 0fill
\noindent{\bf Remark 1} Though ${\mathcal A}^{min}
\subset {\mathcal A}^{com}$, these are not equal in general. For instance, if
${\mathcal B}$ is the ideal of all bounded linear operators, then
${\mathcal B}^{com} = {\mathcal K}$ and ${\mathcal B}^{min} = {\overline{\mathcal F}}$
and ${\mathcal K} \neq {\overline{\mathcal F}}$. However, in many cases
 ${\mathcal A}^{min} = {\mathcal A}^{com}$ and we shall come across some such cases in the paper.
\vskip 8pt plus 0fill
\noindent{\bf Remark 2} We recall the definition of the maximal hull ${\mathcal A}^{max}$ of an operator
ideal $\mathcal A$ and the corresponding procedure $max : {\mathcal A} \to {\mathcal A}^{max}$
\cite[Section 4.9]{P1}. Since  ${\mathcal A}^{min} \subset {\mathcal A}^{com}
\subset {\mathcal A}$ and since $({\mathcal A}^{min})^{max} = {\mathcal A}^{max}$,
we may conclude that $({\mathcal A}^{com})^{max} = {\mathcal A}^{max}$. Also as
${\mathcal A} \subset {\mathcal A}^{max}$, we get
${\mathcal A}^{com} \subset ({\mathcal A}^{max})^{com}$. However, these in general are not equal; for
${\mathcal A} = {\overline{\mathcal F}}$, we get ${\overline{\mathcal F}}\hspace{1mm}^{com} =
{\overline{\mathcal F}}$ (for ${\overline{\mathcal F}} = {\mathcal A}^{min}$) and
$({\overline{\mathcal F}}\hspace{1mm}^{max})^{com} = {\mathcal B}^{com} = {\mathcal K}$ so that
${\overline{\mathcal F}}\hspace{1mm}^{com} \neq ({\overline{\mathcal F}}\hspace{1mm}^{max})^{com}$.
\vskip 8pt plus 0fill

Let $({\mathcal A}, \alpha)$be a quasi-Banach operator ideal. For Banach spaces
$X$ $Y$, an operator $T \in {\mathcal B}(X, Y)$ is said to be in
${\mathcal A}^{reg}(X, Y)$, if $\kappa _Y\circ T \in {\mathcal B}(X, Y^{**})$ and
$\alpha ^{reg}(T) := \alpha(\kappa _Y\circ T)$ \cite[Section 4.5]{P1}.
Now, $({\mathcal A}^{reg}, \alpha ^{reg})$ is a quasi-Banach operator ideal
containing ${\mathcal A}$, the procedure ${\mathcal A} \to {\mathcal A}^{reg}$ is a
hull procedure and ${\mathcal A}^{reg}$ is called the
{\it regular hull} of $\mathcal A$. We now consider the interplay of {\it reg} and
{\it com}.
\begin{proposition}
${\mathcal K}\circ({\mathcal A}^{reg}) = {\mathcal K}\circ{\mathcal A}$. In particular,
$({\mathcal A}^{reg})^{com} = {\mathcal A}^{com}$.
\end{proposition}
\pf
First note that ${\mathcal A} \subset {\mathcal A}^{reg}$, so that
${\mathcal K}\circ{\mathcal A} \subset {\mathcal K}\circ{\mathcal A}^{reg}$. Next, let
$T \in {\mathcal K}\circ{\mathcal A}^{reg}(X, Y)$ for some Banach spaces $X$ and $Y$.
Then there is a Banach space $Z$ and operators $U \in {\mathcal K}(Z, Y)$ and
$S \in {\mathcal A}^{reg}(X, Z)$ such that $T= U\circ S$. Thus ${\kappa}_ZS \in
{\mathcal A}(X, Z^{**})$. As $U$ is compact we get that $U^{**}$ is compact with
$U^{**} (Z^{**}) \subset \kappa _Y(Y)$. Thus we can find $V \in \mathcal{K} (Z^{**}, Y)$ such that $\kappa _y \circ V = U^{**}$. 
Since $\kappa _Y $ is an isometry and $\kappa _Y \circ U = U^{**} \circ \kappa _Z = \kappa _Y \circ V \circ \kappa _Z$, we may conclude that 
$U = V \circ \kappa _Z$. Now it follows that  $T = U\circ S = V  \circ \kappa _Z\circ S \in {\mathcal K} \circ {\mathcal A}(X, Y)$, 
so that ${\mathcal K}\circ({\mathcal A}^{reg}) = {\mathcal K} \circ {\mathcal A}$. Consequently, 
$({\mathcal A}^{reg})^{com} = {\mathcal A}^{com}$. \hfill $\triangle$
\vskip 8pt plus 0fill
\noindent {\bf Remark 1} Trivially, ${\mathcal A}^{com} \subset ({\mathcal A}^{com})^{reg}$.
However, ${\mathcal A}^{com} \neq ({\mathcal A}^{com})^{reg}$. In fact, for ${\mathcal A}
= {\mathcal I}$, the ideal of all integral operators, we have ${\mathcal I}^{com} =
{\mathcal N}$, the ideal of all nuclear operators; and $({\mathcal I}^{com})^{reg} =
{\mathcal N}^{reg}$, but ${\mathcal I}^{com} = {\mathcal N}$ is not regular.

\noindent{\bf Remark 2} A striking difference in the behaviour of $\overline{\mathcal F}$
and $\mathcal K$ in composition is now evident. Indeed, $\overline{\mathcal F} \circ
{\mathcal I} = \overline{\mathcal F}\circ {\mathcal I} \circ \overline{\mathcal F} =
{\mathcal N} = {\mathcal I} \circ \overline{\mathcal F}$, but as $\mathcal N$ is not
regular, ${\mathcal K} \circ
{\mathcal I} = {\mathcal K} \circ {\mathcal I} \circ {\mathcal K} = {\mathcal N} \neq
{\mathcal N}^{reg} = {\mathcal I} \circ {\mathcal K}$ \cite[Theorem 2.1]{Or}.
\vskip 8pt plus 0fill

Let $({\mathcal A}, \alpha)$ be a quasinormed operator ideal. For Banach spaces $X$, $Y$ and operator $T\in {\mathcal B}(X, Y)$ is said to be in ${\mathcal A}^{dual}(X, Y)$ if $T^* \in {\mathcal A}(Y^*, X^*)$ and $\alpha ^{dual}(T) = \alpha(T^*)$ \cite[Section 4.4]{P1}. Now, $({\mathcal A}^{dual}, \alpha ^{dual})$ is a quasinormed operator ideal and is called the {\it dual} of $\mathcal A$. The following observation, which are routine in nature, shall be frequently used in the paper: For an operator ideal $\mathcal{A}$ we have
\begin{enumerate}
\item $\mathcal{A} \cup \mathcal{A}^{dual dual} \subset \mathcal{A}^{reg}$; and
\item $\mathcal{A} \subset  \mathcal{A}^{dual dual}$ if and only if $\mathcal{A}^{dual dual} = \mathcal{A}^{reg}$.
\end{enumerate}
 Next, we study the interplay of the procedures {\it com} and {\it dual}. We note that for any operator ideal ${\mathcal A}$, we have $$({\mathcal A}^{dual})^{com} \subset ({\mathcal A}^{com})^{dual}.$$ But in general, these are not equal; for instance, the ideal of ${\mathcal I}$ of integral operators, ${\mathcal I} = {\mathcal I}^{dual}$ and $({\mathcal I}^{dual})^{com} = {\mathcal I}^{com} = {\mathcal N} \neq {\mathcal N}^{dual} = {({\mathcal I}^{com}})^{dual}$. However,
\vskip 8pt plus 0fill
\begin{proposition}
For an operator ideal $\mathcal A$, with  $\mathcal{A} \subset  \mathcal{A}^{dual dual}$, we have
$$({\mathcal A}^{com})^{dual} = (({\mathcal A}^{dual})^{com})^{reg}.$$
\end{proposition}
\pf
Let $T \in ({\mathcal A}^{com})^{dual}(X, Y)$, for Banach spaces $X$ and $Y$. Then $T^* \in {\mathcal A}^{com}(Y^*, X^*)$. There exists Banach spaces $Z_1$, $Z_2$ and operators $U \in {\mathcal K}(Z_2, X^*)$, $V \in {\mathcal K}(Y^*, Z_1)$, $S \in {\mathcal A} (Z_1, Z_2)$ such that $T^* = U\circ S\circ V$. Since  $\mathcal{A} \subset  \mathcal{A}^{dual dual}$, we get $S^* \in {\mathcal A}^{dual}({Z_2}^*, {Z_1}^*)$ so that $T^{**} = V^*\circ S^*\circ U^* \in ({\mathcal A}^{dual})^{com}(X^{**}, Y^{**})$. It follows that $\kappa _Y\circ T = T^{**}\circ\kappa _X = V^*\circ S^*\circ U^*\circ\kappa _X \in ({\mathcal A}^{dual})^{com}(X, Y^*)$, whence $T \in (({\mathcal A}^{dual})^{com})^{reg}$. Now $({\mathcal A}^{com})^{dual}$ being regular, the result follows. \hfill $\triangle$

\section{Injective and surjective operator ideals - A revisit}

For every quasinormed operator ideal $({\mathcal A}, \alpha)$,
there is a smallest injective operator ideal ${\mathcal A}^{inj}$ containing
$\mathcal A$. For Banach spaces $X$, $Y$ and the natural Alouglu embedding $i_Y : Y \hookrightarrow
l_\infty(B_{Y^*})$, an operator $T \in {\mathcal B}(X, Y)$ is in
${\mathcal A}^{inj}$ if and only if $i_Y\circ T \in {\mathcal A}$ with
$\alpha ^{inj}(T) := \alpha(i_Y\circ T)$. The ideal $({\mathcal A}^{inj}, \alpha ^{inj})$
is a quasinormed operator ideal \cite[Section 4.6]{P1}. The procedure
${\mathcal A} \to {\mathcal A}^{inj}$ is a hull procedure and ${\mathcal A}^{inj}$ is
called the injective hull of $\mathcal A$. The operator ideal $\mathcal A$ is said to
be injective if ${\mathcal A}^{inj} = {\mathcal A}$.

First of all we give an interesting and useful characterization of injective operator ideals.

\begin{definition}
An operator ideal ${\mathcal A}$ is said to have the restricted range property ({\it RRP},
for short), if for arbitrary Banach spaces $X$, $Y$ and $T \in {\mathcal A}(X, Y)$, we have $\hat{T} \in {\mathcal A}(X, {\overline{R(X)}})$.
Here $R(X)$ is the range of $T$ in $Y$ and $\hat {T}(x) = T(x)$ for all $x \in X$.
\end{definition}
\begin{lemma}
An operator ideal ${\mathcal A}$ is injective if and only if it has the restricted range
property.
\end{lemma}
\pf
It follows from \cite[Proposition 8.5.4]{P1} that an injective operator ideal has the {\it RRP}.

Conversely, let ${\mathcal A}$ have the {\it RRP}. Let $T \in {\mathcal A}^{inj}(X, Y)$
for some Banach spaces. Then $i_Y\circ T \in {\mathcal A}(X, l_\infty(B_{Y^*}))$.
Since ${\mathcal A}$ has the {\it RRP}, we have $T_1 \in {\mathcal A}
(X, \overline{i_Y\circ T(X)})$, where $T_1(x) = i_Y\circ T(x)$, $x\in X$. As $i_Y$
is an isometry and as $\overline{T(X)} \subset Y$ we may obtain an isometry
$I_0: {\mathcal B}(\overline{i_y\circ T(X)}, Y)$ given by $I_0(i_Y\circ T(x)) = Tx$,
$x\in X$. Now $I_0\circ T_1(x) = I_0(i_y\circ T(x)) = T(x)$ for all $x\in X$.
In other words, $T = I_0\circ T_1 \in {\mathcal A}(X, Y)$ so that $\mathcal A$ is injective.
   \hfill $\triangle$

\begin{theorem}
If ${\mathcal A}_1$ and ${\mathcal A}_2$ are injective operator ideals, then so is
${\mathcal A}_1\circ{\mathcal A}_2$.
\end{theorem}

\pf Let $T \in ({\mathcal A}_1\circ{\mathcal A}_2)^{inj}
(X, Y)$ and let $I: Y \to Z$ be an isometry for some Banach spaces $X$, $Y$ and $Z$.
Then $I\circ T \in ({\mathcal A}_1\circ {\mathcal A}_2)(X,Z)$. Thus there exists a
Banach space $W$ and operators $U \in {\mathcal A}_1(W, Z)$ and $V \in {\mathcal A}_2
(X, W)$ such that $I\circ T = U\circ V$. Since ${\mathcal A}_2$ being injective has
the {\it RRP}, by using the ideal property of ${\mathcal A}_1$ we may assume that
$W = \overline{V(X)}$. Next, we set $U_0: V(X) \to Y$ by $U_0(Vx) = Tx$, for all
$X \in X$. If $Vx_0 = 0$ for some $x_0 \in X$, then  $I\circ T(x_0) = U\circ V(x_0) = 0$. Since
$I$ is an isometry we conclude that $Tx = 0$. Thus $U_0$ is a well defined linear map.
Also $$\|U_0(Vx)\| = \|Tx \| = \|I\circ (x)\| = \|U\circ V(x)\| \leq \|U\|\|Vx\|,$$
for all $x \in X$, so that $U_0 \in {\mathcal B}(V(X), Y)$. Thus we can extend $U_0$
to $U_1 \in {\mathcal B}(W, Y)$. Now
$$I\circ U_1(Vx) = I\circ T(x) = U(Vx), for \ all \ x \in X$$
so that $U = I\circ U_1 \in {\mathcal A}_1(W, Z)$. As ${\mathcal A}_1$ is injective,
we conclude that $U_1 \in {\mathcal A}_1(W, Z)$. It follows that $T = U_1\circ V \in
({\mathcal A}_1\circ{\mathcal A}_2)(X, Y)$. In other words ${\mathcal A}_1\circ{\mathcal A}_2$
is injective.
   \hfill $\triangle$
\vskip 8pt plus 0fill
\noindent{\bf Remark} If ${\mathcal A}$ is an injective operator ideal then so are ${\mathcal K}\circ{\mathcal A}$ and ${\mathcal A}\circ{\mathcal K}$. In particular, ${\mathcal A}^{com}$ is injective whenever ${\mathcal A}$ is so.

\vskip 8pt plus 0fill
For every quasinormed operator ideal $({\mathcal A}, \alpha)$,
there is a smallest operator ideal ${\mathcal A}^{sur}$  containing
$\mathcal A$. For Banach spaces $X$ and $Y$, an operator $T \in {\mathcal B}(X, Y)$ is in
${\mathcal A}^{sur}$ if and only if $T\circ q_X \in {\mathcal A}$ with
$\alpha ^{sur}(T) := \alpha(T\circ q_X)$. The ideal $({\mathcal A}^{sur}, \alpha ^{sur})$
is a quasinormed operator ideal \cite[Section 4.7]{P1}. The procedure
${\mathcal A} \to {\mathcal A}^{sur}$ is a hull procedure and ${\mathcal A}^{sur}$
is called the surjective hull of $\mathcal A$. The operator ideal $\mathcal A$
is called surjective if ${\mathcal A}^{sur} = {\mathcal A}$.
We can dualise the above ideas for surjective operator ideals as follows.

To prove this theorem we shall use the following characterization of surjective
operator ideals.
\begin{definition}
An operator ideal ${\mathcal A}$ is said to have the quotiented domain property
(QDP, for short) if for arbitrary Banach spaces $X$, $Y$ and operator
$T \in {\mathcal A}(X, Y)$, we have $T_0 \in {\mathcal A}(X/Ker \ T, Y)$is the kernel of T
and $T_0(X + Ker \ T) = Tx$, for all $x \in x$.
\end{definition}

\begin{lemma}An operator ideal
${\mathcal A}$ is surjective if and only if it has quotiented domain property.
\end{lemma}
\pf
It follows from \cite[proposition 8.5.4]{P1} that a surjective operator ideal has the {\it QDP}.

Conversely, assume that $\mathcal A$ has the {\it QDP} and for a pair of Banach spaces
$X$ and $Y$, let $T \in {\mathcal A}^{sur}
(X, Y)$. Then $T\circ q_X \in {\mathcal A}(l_1(B_X), Y)$. Let $Z = l_1(B_X)/Ker(T\circ q_X)$.
Let $T_0: Z \to Y$ be the map corresponding to $T\circ q_X$ and let $q: l_1(B_X) \to Z$ be
the natural quotient map. Then $T\circ q_X = T_0\circ q$, and we have $T_0 \in {\mathcal A}(Z, Y)$ for
$\mathcal A$ has the {\it QDP}. For each $x \in X$ we set $V_0(x) = q(\alpha)$, where
$x = q_X(\alpha)$. Now if $q_X(\alpha) = 0$, then $Tq_X(\alpha) = 0$ so that $\alpha
\in Ker(T\circ q_X)$. Thus $V_0(0) = q (0) = 0$. Therefore, $V_0: X \to Z$ is a well
defined linear operator. Also, for any $x \in X$ we have
$$\|V_0(x)\| = {\textrm inf}\{\|q_X(\alpha)\|: x = q_X(\alpha)\} \leq {\textrm inf}\{\|\alpha\|: x = q_X(\alpha)\} = \|x\|.$$
Thus $V_0 \in {\mathcal B}(X, Z)$. Now $T_0\circ V_0(x) = T_0\circ q(\alpha) = T\circ q_X(\alpha) =
T(x)$ for all $x \in X$, so that $T = T_0\circ V_0 \in {\mathcal A}(X, Y)$. Hence
$\mathcal A$ is surjective.
   \hfill $\triangle$

\begin{theorem}
If ${\mathcal A}_1$ and ${\mathcal A}_2$ are surjective operator ideals, then so is
${\mathcal A}_1\circ{\mathcal A}_2$.
\end{theorem}

\pf Let $X$, $Y$ be Banach spaces and $T \in
({\mathcal A}_1\circ{\mathcal A}_2)^{sur}(X, Y)$. Then $T\circ q_X \in ({\mathcal A}_1\circ{\mathcal A}_2)(l_1(B_X), Y)$. Thus there
is a Banach space $Z$ and operators $U \in {\mathcal A}_1(Z, Y)$ and $V \in
{\mathcal A}_2(l_1(B_X), Z)$ such that $T\circ q_X = U\circ V$. Let $q_U: Z \to Z/Ker U$
be the natural quotient map and $U_0: Z/Ker \ U \to Y$ be the map corresponding to $U$.
Then, $U = U_0\circ q_U$. Since $\mathcal A$ has the {\it QDP}, we have
$U_0 \in {\mathcal A}(Z/Ker \ U, Y)$. For $x \in X$ we set $V_0(x) = q_UV(\alpha)$, where
$x = q_X(\alpha)$. Then, proceeding as in Lemma 3.5, we may obtain $V_0 \in
{\mathcal B}(X, Z/Ker U)$ and that $V_0\circ q_X = q_U\circ V \in {\mathcal A}_2
(l_1(B_X), Z/Ker U)$. Since ${\mathcal A}_2$ is surjective, we conclude further
that $V_0 \in {\mathcal A}_2(X, Z/Ker U)$. Now for all $x \in X$
$$U_0\circ V_0(x) = U_0(q_U\circ V(\alpha)) = U\circ V(\alpha) = T\circ q_X(\alpha) = T(x),$$
so that $T = U_0\circ V_0 \in {\mathcal A}_1\circ{\mathcal A}_2(X, Y)$.
Hence ${\mathcal A}_1\circ{\mathcal A}_2$ is surjective.
   \hfill $\triangle$
\vskip 8pt plus 0fill
\noindent{\bf Remark} Let $\mathcal A$ be a surjective operator ideal, then so
are ${\mathcal K}\circ{\mathcal A}$ and ${\mathcal A}\circ{\mathcal K}$. In particular,
${\mathcal A}^{com}$ is surjective if $\mathcal A$ is so.

\section{ Extension and lifting properties for operator ideals}
Let ${\mathcal A}_1$ and ${\mathcal A}_2$ be two operator ideals. Then by Theorem 3.3,
${{\mathcal A}_1}^{inj}\circ{{\mathcal A}_2}^{inj}$ is an injective operator ideal
containing $({\mathcal A}_1 \circ{\mathcal A}_2)^{inj}$. Also it follows from the definition
of injectivity, that ${{\mathcal A}_1}^{inj} \circ{\mathcal A}_2 \subset ({\mathcal A}_1
\circ{\mathcal A}_2)^{inj}$. Thus
$${{\mathcal A}_1}^{inj} \circ{\mathcal A}_2 \subset ({\mathcal A}_1
\circ{\mathcal A}_2)^{inj} \subset {\mathcal A}_1^{inj}
\circ{\mathcal A}_2^{inj}.$$
\noindent In particular, when ${\mathcal A}_2$ is injective, we have
$$({\mathcal A}_1
\circ{\mathcal A}_2)^{inj} = {{\mathcal A}_1}^{inj} \circ{\mathcal A}_2.$$
\noindent Next, we propose a condition on ${\mathcal A}_1$ for which
$({\mathcal A}_1
\circ{\mathcal A}_2)^{inj} = {{\mathcal A}_1}^{inj} \circ{{\mathcal A}_2}^{inj}$.

\begin{definition}
An operator ideal $({\mathcal A}, \alpha)$ is said to have the $l_\infty$-extension property
if for a Banach spaces $X$, a set $\Gamma$ and an operator $T \in {\mathcal A}
(X, l_\infty(\Gamma))$, there is an operator ${\tilde{T}} \in {\mathcal A}(l_\infty(B_{X^*}), l_\infty(\Gamma))$
such that $T = {\tilde{T}}\circ i_X$, with $\alpha({\tilde{T}}) = \alpha(T)$, where
$i_X: X \hookrightarrow l_\infty (B_{X^*})$ is the Alouglu embedding.
\end{definition}

Since $l_\infty(\Gamma)$ spaces are injective, the ideals $\mathcal B$ and $\mathcal K$
of all bounded and all compact operators respectively, have the $l_\infty$-extension property.
Since $l_\infty(\Gamma)$ has the approximation property, ${\overline{\mathcal F}}
(X, l_\infty(\Gamma)) = {\mathcal K}(X, l_\infty(\Gamma))$ for all Banach spaces $X$. It
follows that $\overline{\mathcal F}$ has the $l_\infty$-extension property.
The ideal ${\mathcal I}_p$ of $p$-integral operators for $1 \leq p \leq \infty$ has this
property. In fact,
these operator ideals enjoy a much stronger extension property \cite[Proposition 6.12]{DJT}.
As for the operator ideals $\Pi _p$, $\Pi _p(X, Y)$ coincides with ${\mathcal I}_p(X, Y)$
if $Y$ is an $l_\infty$-space. It follows that $\Pi _p$ also enjoys the $l_\infty$-
extension property.

\begin{theorem}
Let ${\mathcal A}_1$ and ${\mathcal A}_2$ be two operator ideals. If
${\mathcal A}_1$ has the $l_\infty$-extension property, or if ${\mathcal A}_2$
is injective, then
$$({\mathcal A}_1
\circ{\mathcal A}_2)^{inj} = {{\mathcal A}_1}^{inj} \circ{{\mathcal A}_2}^{inj}.$$
\end{theorem}

\pf Let $X$ and $Y$ be Banach spaces and let $T \in {{\mathcal A}_1}^{inj}
\circ{{\mathcal A}_2}^{inj}(X, Y)$. Then there is a Banach space $Z$ and operators
$U \in {{\mathcal A}_1}^{inj}(Z, Y)$ and $V \in {{\mathcal A}_2}^{inj}(X, Z)$
such that $T = U \circ V$. If ${\mathcal A}_1$ has the $l_\infty$-extension property,
then there is an $\tilde{U} \in {\mathcal A}_1(l_\infty(B_{Z^*}), l_\infty(B_{Y^*}))$
such that $i_Y \circ U = \tilde{U} \circ i_Z$. Thus we have
$$i_Y \circ T = i_Y \circ U \circ V = \tilde{U} \circ i_Z\circ V =
\tilde{U}\circ\tilde{V},$$
where $\tilde{V} = i_Z\circ V \in {\mathcal A}_2(X, l_\infty(B_{Z^*})$. It
follows that $T \in {{\mathcal A}_1}^{inj} \circ{{\mathcal A}_2}^{inj}(X, Y)$.
The remaining part of the proof follows from the above discussion.
   \hfill $\triangle$

\begin{proposition}
If ${\mathcal A}_1$ and ${\mathcal A}_2$ are two operators ideals both having
the $l_\infty$-extension property, then so has ${\mathcal A}_1\circ{\mathcal A}_2$.
\end{proposition}

\pf
Let $X$ be a Banach space, $\Gamma$ a set and consider the operator $T \in
{\mathcal A}_1\circ{\mathcal A}_2(X, l_\infty(\Gamma))$. Then there is are
a Banach space $Y$ and operators $U \in {\mathcal A}_1(Y, l_\infty(\Gamma))$
and $V \in {\mathcal A}_2(X, Y)$ such that $T = U \circ V$. Since
both ${\mathcal A}_1$ and ${\mathcal A}_2$ have the $l_\infty$-extension property,
$U$ has an extension ${\tilde{U}} \in {\mathcal A}_1(l_\infty(B_{Y^*}), l_\infty(\Gamma))$
and $i_Y\circ V$ has an extension
${\tilde{V}} \in {\mathcal A}_2(l_\infty(B_{X^*}), l_\infty(B_{Y^*})$.
Thus, ${\tilde{U}}\circ{\tilde{V}}\circ i_X =
{\tilde{U}}\circ i_Y\circ V = U\circ V$. Hence, ${\tilde{U}}\circ{\tilde{V}} \in
{\mathcal A}_1\circ{\mathcal A}_2$ extends $T$.
   \hfill $\triangle$

Let ${\mathcal A}_1$ and ${\mathcal A}_2$ be two operator ideals. It follows by
Theorem 3.6 that
$${\mathcal A}_1\circ{\mathcal A}_2^{sur} \subset ({\mathcal A}_1\circ{\mathcal A}_2)^{sur}
\subset {\mathcal A}_1^{sur}\circ{\mathcal A}_2^{sur}.$$
In particular, if ${\mathcal A}_1$ is surjective then $({\mathcal A}_1\circ{\mathcal A}_2)^{sur}
= {\mathcal A}_1\circ{\mathcal A}_2^{sur}$. We propose, once more, a condition on
${\mathcal A}_2$ such that $({\mathcal A}_1\circ{\mathcal A}_2)^{sur}
= {\mathcal A}_1^{sur}\circ{\mathcal A}_2^{sur}$.

\begin{definition}
An operator ideal $({\mathcal A}, \alpha)$ is said to have the $l_1$-lifting property
if for any set $\Gamma$, a Banach space $Y$ and an operator $T \in {\mathcal A}
(l_1(\Gamma), Y)$, there is an ${\hat{T}} \in (l_1(\Gamma), l_1(B_Y))$
such that $q_Y\circ{\hat{T}} = \kappa _Y\circ T$, with $\alpha({\hat{T}}) = \alpha(T)$,
where $\kappa _Y:Y \hookrightarrow  Y^{**}$ is the canonical embedding and
$q_Y: l_1(B_Y) \to Y$ is the natural quotient map.
\end{definition}

The lifting property of $l_1(\Gamma)$ ensures that both the operator ideals $\mathcal B$
and $\mathcal K$ have the $l_1$-lifting property. Since the $l_1^*(\Gamma)$ has the
approximation property, $\overline{\mathcal F}(l_1(\Gamma), Y) = {\mathcal K}(l_1(\Gamma), Y)$
for all Banach spaces $Y$, it follows that $\overline{\mathcal F}$ has the
$l_1$-lifting property. Using the $l_\infty$-extension property of
$\Pi _p$ it is easy to verify that $\Pi _p^{dual}$ also has the $l_1$-lifting property.
Dualising the proof of Theorem 4.2, we can prove the following.

\begin{theorem}
Let ${\mathcal A}_1$ and ${\mathcal A}_2$ be two operator ideals. If
${\mathcal A}_1$ is surjective, or if ${\mathcal A}_2$ has the $l_1$-lifting property,
then
$$({\mathcal A}_1
\circ{\mathcal A}_2)^{sur} = {{\mathcal A}_1}^{sur} \circ{{\mathcal A}_2}^{sur}.$$
\end{theorem}

\begin{proposition}
If ${\mathcal A}_1$ and ${\mathcal A}_2$ are two operators ideals both having
the $l_1$-lifting property, then so has ${\mathcal A}_1\circ{\mathcal A}_2$.
\end{proposition}
\pf
Let $Y$ be a Banach space, $\Gamma$ a set and $T \in {\mathcal A}_1\circ
{\mathcal A}_2 (l_1(\Gamma), Y)$. Then, there is a Banach space $Z$ and operators
$U \in {\mathcal A}_1(X, Y)$ and $V \in {\mathcal A}_2(l_1(\Gamma), Z)$ such that
$T = U \circ V$. Since, ${\mathcal A}_2$ has the $l_1$-lifting property, there is
an operator $\hat{V} \in {\mathcal A}_2(l_1(\Gamma), l_1(B_Z))$ such that
$V = q_Z \circ \hat{V}$.
Next, $U \circ q_Z \in {\mathcal A}_1(l_1(B_Z), Y)$ and ${\mathcal A}_1$ has the
$l_1$-lifting property. Thus there exists an operator
$\hat{U} \in {\mathcal A}_1(l_1(B_Z), l_1(B_Y))$ such that
$U \circ q_Z = q_Y \circ \hat{U}$. Thus $$T = U \circ V = U \circ q_Z \circ \hat{V}
= q_Y \circ \hat{U} \circ \hat{V}.$$
Now, $\hat{U} \circ \hat{V} \in {\mathcal A}_1\circ{\mathcal A}_2(l_1(\Gamma), l_1(B_Y))$,
so that we may conclude that ${\mathcal A}_1\circ{\mathcal A}_2$ has the $l_1$-lifting
property.
   \hfill $\triangle$
\vskip 8pt plus 0fill
Now we investigate the duality relationship between the
$l_\infty$-extension and the $l_1$-lifting properties. To begin with we prove two lemmas.

\begin{lemma}
Let $X$ be a Banach space, $\Gamma$ a set and $T\in {\mathcal A}(X, l_\infty(\Gamma))$.
Then the following are equivalent:

{\leftskip=6mm
\item{(a)}
For any Banach space $Y$, any isometry $I: X \to Y$ and $\epsilon > 0$, there is a
$\tilde{T} \in {\mathcal A}(Y, l_\infty(\Gamma))$ with $\alpha(\tilde{T}) \leq
(1 + \epsilon)\alpha(T)$ such that $T = \tilde{T} \circ I$.
\item{(b)}
For some set $\Omega$, and isometry $i_0: X \to l_\infty(\Omega)$ and $\epsilon > 0$,
there is a $T_0 \in {\mathcal A}(l_\infty(\Omega), l_\infty(\Gamma))$ with
$\alpha(\tilde{T}) \leq (1 + \epsilon) \alpha(T)$ such that $T = T_0 \circ i_0$. \par}
\end{lemma}

\pf
It suffices to show that $(b) \Rightarrow (a)$. Let $I: X \to Y$ be any isometry.
Consider the canonical embedding $i_X: X \hookrightarrow l_\infty(B_{X^*})$. Since
$l_\infty(B_{X^*})$ has the extension property, there is a bounded operator $\tilde{i_X}: Y
\to l_\infty(B_{X^*})$ with $\|\tilde{i_X}\| \leq (1 + \epsilon /4)$ such that $i_X = \tilde{i_X}\circ I$.
Now by (b), $T$ has an extension $T_0: l_\infty(B_{X^*}) \to l_\infty(\Gamma)$ with $\alpha(T_0) < (1 + \epsilon /4)
\alpha(T)$ such that  $T = T_0 \circ i_X = \tilde{T}\circ I$ where $\tilde{T} = T_0 \circ
\tilde{i_X} \in {\mathcal A}(Y, l_\infty(\Gamma))$.
   \hfill $\triangle$

\begin{lemma}
Let $Y$ be a Banach space, $\Gamma$ a set and $T \in {\mathcal A}(l_1(\Gamma), Y)$.
Then the following are equivalent: \par

{\leftskip=6mm
\item{(a)}
For any Banach space $X$, metric quotient $Q: X \to Y$ and $\epsilon > 0$, there is a
$\hat{T} \in {\mathcal A}((l_1(\Gamma), Y)$ with $\alpha(\hat{T}) \leq (1 + \epsilon)\alpha(T)$
such that $T = Q\circ\hat{T}$
\item{(b)}
For some set $\Omega$, a metric quotient $q_0: l_1(\Omega) \to Y$ and $\epsilon > 0$,
there is a $T_0 \in {\mathcal A}(l_1(\Gamma), l_1(\Omega))$ with $\alpha(T_0) \leq
(1 + \epsilon)\alpha(T)$ such that $T = q_0 \circ T_o$. \par}
\end{lemma}

\pf
It suffices to prove that $(b) \Rightarrow (a)$. Let $Q: X \to Y$ be any metric quotient.
Consider the natural metric quotient $q_X: l_1(B_X) \to X$.
Since $l_1(B_X)$ has the lifting property, there is a bounded operator $\hat{q_X}:
l_1(B_X) \to Y$ with $\|\hat{q_X}\| < 1 + \epsilon /4$ such that $\hat{q_X} = Q \circ q_X$.
Now by (b), $T$ has a lifting $T_0: {\mathcal A}(l_1(\Gamma), l_1(B_X))$ with
$\alpha(T_0) < (1 + \epsilon /4) \alpha (T)$ such that
$T = \hat{q_X} \circ T_0 = Q \circ  \hat{T}$ where $\hat{T} = q_X \circ T_0 \in
{\mathcal A}(l_1(\Gamma), X))$.
   \hfill $\triangle$

\vskip 8pt plus 0fill
\noindent{\bf Remark} Since ${\mathcal A}(X, l_\infty(\Gamma)) =
{\mathcal A}^{inj}(X, l_\infty(\Gamma)) = {\mathcal A}^{reg}(X, l_\infty(\Gamma))$, we conclude that
${\mathcal A}$ has the $l_\infty$-extension property if and only if ${\mathcal A}^{inj}$ has the
$l_\infty$-extension property if and only if ${\mathcal A}^{reg}$has the
$l_\infty$-extension property. Dually, since ${\mathcal A}(l_1(\Gamma), Y)
= {\mathcal A}^{sur}(l_1(\Gamma), Y)$, we may further conclude that $\mathcal A$ has
the $l_1$-lifting property if and only of so does ${\mathcal A}^{sur}$. However, we
are not aware of the dual situation for the case of the regular hull.
\vskip 8pt plus 0fill

\begin{theorem}
If $\mathcal A$ has the $l_\infty$-extension property, then ${\mathcal A}^{dual}$
has the $l_1$-lifting property. If ${\mathcal A} \subset {\mathcal A}^{dual\hspace{2mm}dual}$
(or equivalently, ${\mathcal A}^{dual\hspace{2mm}dual} = {\mathcal A}^{reg})$), then the
converse also holds
\end{theorem}

\pf
Let $Y$ be a Banach space, $\Gamma$ a set and $T \in {\mathcal A}^{dual}
(l_1(\Gamma), Y)$. Then $T^* \in {\mathcal A} (Y^*, l_\infty(\Gamma))$. Since
$\mathcal A$ has the $l_\infty$-extension property, there is an
$S_0 \in {\mathcal A}(l_\infty(B_{Y^*}), l_\infty(\Gamma))$
such that $T^* = S_0 \circ q_Y^*$. Consider he adjoint $S_0^*: l_1(\Gamma)^{**}
\to l_1(B_{Y^*})^*$ of $S_0$, and put
$\hat{T} = P_{l_1(B_X)}\circ S_0^*\circ \kappa _{l_1(\Gamma)}$. Then $\hat{T^*} = S_0$,
so that $\hat{T} \in {\mathcal A}^{dual}(l_1(\Gamma), l_1(B_{Y^*}))$, and we have
$T = q_Y\circ\hat{T}$. Hence ${\mathcal A}^{dual}$ has the $l_1$-lifting property.

Conversely, let ${\mathcal A} \subset {\mathcal A}^{dual\hspace{2mm}dual}$ and let
${\mathcal A}^{dual}$ have the $l_1$-lifting property. For a Banach space $X$ and
a set $\Gamma$, let $T \in {\mathcal A}(X, l_\infty(\Gamma)) \subset
{\mathcal A}^{dual}(X, l_\infty(\Gamma))$. Then $T^* \in {\mathcal A}^{dual}
(l_1(\Gamma)^{**}, X^*)$, so that $T^*\circ\kappa _{l_1(\Gamma)} \in {\mathcal A}^{dual}
(l_1(\Gamma), X^*)$. Since ${\mathcal A}^{dual}$ has the $l_1$-lifting property,
there is a $S_0 \in {\mathcal A}^{dual}(l_1(\Gamma), l_1(B_{X^*}))$ such that
$T^*\circ \kappa _{l_1(\Gamma)} = q_X \circ S$. Thus
$$T = \kappa _{l_1(\Gamma)}^* \circ\kappa _{l_\infty(\Gamma)}\circ T =
\kappa _{l_1(\Gamma)}^* \circ T^{**}\circ\kappa _X = S_0^*\circ q_{X^*}^*\circ\kappa _X
= S_0^*\circ i_X.$$
Thus, $S_0^* \in {\mathcal A}({l_\infty(B_{X^*})}, l_\infty(\Gamma))$. Hence
$\mathcal A$ has the $l_\infty$-extension property.   \hfill $\triangle$

\begin{theorem}
Let ${\mathcal A} \subset {\mathcal A}^{dual\hspace{2mm}dual}$
(or equivalently,  ${\mathcal A}^{dual\hspace{2mm}dual} = {\mathcal A}^{reg})$. Then,
 $\mathcal A$ has the $l_1$-lifting property if and only if
${\mathcal A}^{dual}$ has the $l_\infty$-extension property.
\end{theorem}

\pf
First assume that $\mathcal A$ has the $l_1$-lifting property. For a Banach space $X$
and a set $\Gamma$, let $T \in {\mathcal A}^{dual}(X, l_\infty(\Gamma))$. Then
$T^* \in {\mathcal A}(l_1(\Gamma)^{**}, X^*)$ so that $T^*\circ \kappa _{l_1(\Gamma)}
\in {\mathcal A}(l_1(\Gamma), X^*)$. Since $\mathcal A$ has the $l_1$-lifting property, there is an
$S_0 \in{\mathcal A}(l_1(\Gamma), l_1(B_{X^*}))$ such that $T^*\circ \kappa _{l_1(\Gamma)}
= q_{X^*}\circ S_0$. Then
$$T = \kappa _{l_(\Gamma)}^*\circ \kappa _{l_\infty(\Gamma)}\circ T =
\kappa _{l_1(\Gamma)}^*\circ T^{**}\circ \kappa _X = S_0^*\circ q_{X^*}^*\circ \kappa _X
= S_0^*\circ i_X,$$
for $\kappa _{l_1(\Gamma)}^* = P_{l_\infty(\Gamma)}$, $P_{l_\infty(\Gamma)}\circ
\kappa _{l_\infty(\Gamma)} = I_{l_\infty(\Gamma)}$ and $q_{X^*}^*\circ \kappa _X =i_X$.
Put $S_0^* = \hat{T}$. Since $S_0 \in {\mathcal A}(l_1(\Gamma), l_\infty(B_{X^*}))
\subset {\mathcal A}^{dual \hspace{2mm}dual}(l_1(\Gamma), l_1(B_{X^*}))$, we get
$\hat{T} \in {\mathcal A}^{dual}(l_\infty(B_{X^*}), l_\infty(\Gamma))$. Thus
${\mathcal A}^{dual}$ has the $l_\infty$-extension property.

Conversely, assume that ${\mathcal A}^{dual}$ has the $l_\infty$-extension property.
Then by Theorem 4.9, ${\mathcal A}^{dual \hspace{2mm}dual}$ has the $l_1$-lifting
property. Let $T \in {\mathcal A}(l_1(\Gamma), l_1(B_Y)) \subset
{\mathcal A}^{dual\hspace{2mm}dual}(l_1(\Gamma), l_1(B_Y))$, such that
$T = q_Y\circ \hat{T}$. Now ${\hat{T}}^{**} \in {\mathcal A}(l_1(\Gamma)^{**}, l_1(B_Y)^{**})$
so that $\hat{T} = P_{l_1(B_X)}\circ\kappa _{l_1(B_X)}\circ \hat{T} =
P_{l_1(B_X)} \circ {\hat{T}}^{**} \circ \kappa _{l_1(\Gamma)} \in
{\mathcal A}(l_1(\Gamma), l_1(B_Y))$. Therefore, $\mathcal A$ has the $l_1$-lifting property.
   \hfill $\triangle$

\section{Compact kernels of injective and surjective operator ideals}

In this section we record the interplay of the kernel procedure {\it `com'}
with the hull procedures {\it `inj'} and {\it `sur'}.  The last two sections of preparation
leads us to several observations.

\begin{proposition} Let $\mathcal A$ be an operator ideal. The following are in order \par
{\leftskip=6mm
\item{(1)} Since $\mathcal K$ is injective, we have
$$({\mathcal A}\circ{\mathcal K})^{inj} = {\mathcal A}^{inj}\circ{\mathcal K}.$$
[Theorem 3.3]

\item{(1$^\prime$)} If ${\mathcal A}$ has the $l_\infty$-extension property, then in view of
 ${\mathcal K} = \overline{\mathcal F}\hspace{1mm}^{inj}$, we have
$${\mathcal A}^{inj}\circ{\mathcal K} = ({\mathcal A}\circ\overline{\mathcal F})^{inj}.$$
[Theorem 4.2]
\item{(2)} Since ${\mathcal K}$ is injective and also has the $l_\infty$-extension
property, we have  $${\mathcal K}\circ{\mathcal A}^{inj} = ({\mathcal K}\circ{\mathcal A})^{inj}.$$
[Theorem 4.2]
\item{(2$^\prime$)} In particular, since ${\mathcal K} = \overline{\mathcal F}\hspace{1mm}^{inj}$ and
$\overline{\mathcal F}$ has the $l_\infty$-extension property, we have
$${\mathcal K}\circ{\mathcal A}^{inj} = \overline{\mathcal F}\hspace{1mm}^{inj}\circ{\mathcal A}^{inj}
= (\overline{\mathcal F}\circ{\mathcal A})^{inj}.$$

\item{(3)} It follows from 1 and 2 above, that $$({\mathcal A}^{com})^{inj}
= ({\mathcal A}^{inj})^{com}.$$

\item{(3$^\prime$)} If ${\mathcal A}$ has the $l_\infty$-extension property, then
by 1$^\prime$ and 2$^\prime$ above, we have
$$({\mathcal A}^{inj})^{com}
= {\mathcal K}\circ({\mathcal A}\circ\overline{\mathcal F})^{inj} =
(\overline{\mathcal F}\circ{\mathcal A}\circ\overline{\mathcal F})^{inj} =
({\mathcal A}^{min})^{inj}.$$}\par
\end{proposition}
\vskip 8pt plus 0fill
\noindent{\bf Remarks 1.} Since $\mathcal{I}$ has the $l_{\infty}$-extension property for $1 \le p < \infty$, Proposition 19.2.16 in \cite{P1} follows from Proposition 5.1(3$^\prime$) above.
Also, note that for the operator ideal ${\mathcal K}_p$ defined in \cite{SK1}, by Proposition 5.1(1) above, we have ${\mathcal K}_p^{dual}={\Pi _p}\circ{\mathcal K}={\mathcal I}_p^{inj}\circ{\mathcal K}=({\mathcal I}_p\circ{\mathcal K})^{inj}=
{\mathcal N}_p^{inj}={\mathcal Q}{\mathcal N}_p$, where ${\mathcal Q}{\mathcal N}_p$ is the operator ideal of $p$-quasi nuclear operators.

{\bf 2.} Since $\Pi _p$ is a left accessible injective operator ideal with the $l_\infty$-extension property
$\Pi _p^{dual}$ is surjective, and by Theorem 4.9 it also has the
$l_1$-lifting property.
Thus the results $({\Pi _p}^{min})^{inj} = {\mathcal K}_p^{dual} = {\Pi _p} \circ {\mathcal K}$ and
$(({\Pi _p}^{dual})^{min})^{sur} = {\mathcal K}_p^{dual \hspace{1mm}dual} =
{\mathcal K} \circ {\Pi _p}^{dual}$ for $1 \leq p \leq \infty$,
obtained in \cite{SK2} provide us specific examples
of situations as clarified by the above proposition. Also note that in the limiting case
${\mathcal B}\circ{\mathcal K} = {\mathcal K}\circ{\mathcal B} = {\mathcal K}$ holds
trivially and also we have $({\mathcal B}^{min})^{inj} = \overline{\mathcal F}^{inj} =
{\mathcal K} = ({\mathcal B}^{min})^{sur} = \overline{\mathcal F}^{sur}$.
\vskip 8pt plus 0fill
\begin{proposition} Let $\mathcal A$ be an operator ideal. The following are in order \par
{\leftskip=6mm
\item{(1)} Since $\mathcal K$ is surjective, we have
$$({\mathcal K}\circ{\mathcal A})^{sur} = {\mathcal K}\circ{\mathcal A}^{sur}.$$
[Theorem 3.6]
\item{(1$^\prime$)} If ${\mathcal A}$ has the $l_1$-lifting property, then we have
$${\mathcal K}\circ{\mathcal A}^{sur} = (\overline{\mathcal F}\circ{\mathcal A})^{sur}.$$
[Theorem 4.5]
\item{(2)} Since ${\mathcal K}$ is surjective and also has the $l_1$-lifting
property, we have  $${\mathcal A}^{sur}\circ{\mathcal K} = ({\mathcal A}\circ{\mathcal K})^{sur}.$$
[Theorem 4.5]
\item{(2$^\prime$)} In particular, since ${\mathcal K} = \overline{\mathcal F}\hspace{1mm}^{sur}$ and
$\overline{\mathcal F}$ has the $l_1$-lifting property, we have
$${\mathcal A}^{sur}\circ{\mathcal K} = {\mathcal A}^{sur}\circ\overline{\mathcal F}\hspace{1mm}^{sur}
= ({\mathcal A}\circ\overline{\mathcal F})^{sur}.$$

\item{(3)} It follows from 1 and 2 above, that $$({\mathcal A}^{com})^{sur}
= ({\mathcal A}^{sur})^{com}.$$

\item{(3$^\prime$)} If ${\mathcal A}$ has the $l_1$-lifting property, then
by 1$^\prime$ and 2$^\prime$ above, we have
$$({\mathcal A}^{sur})^{com}
= (\overline{\mathcal F}\circ{\mathcal A})^{sur}\circ {\mathcal K} =
(\overline{\mathcal F}\circ{\mathcal A}\circ\overline{\mathcal F})^{sur} =
({\mathcal A}^{min})^{sur}.$$}\par
\end{proposition}
Let us rewrite the above results for accessible operator ideals.
\begin{corollary}
Let $\mathcal A$ be an operator ideal. Then \par
{\leftskip=6mm
\item{(a$_1$)}If $\mathcal A$ is left accessible, then
$$ ({\mathcal A}^{com})^{inj} = ({\mathcal A}^{min})^{inj} =
{\mathcal K}\circ{\mathcal A}^{inj}.$$
\item{(a$_2$)}If ${\mathcal A}$ is right accessible and also has the $l_\infty$-extension
property, then we have
$$({\mathcal A}^{com})^{inj} = ({\mathcal A}^{min})^{inj} =
{\mathcal A}^{inj}\circ{\mathcal K}.$$

\item{(b$_1$)}If $\mathcal A$ is right accessible, then
$$ ({\mathcal A}^{com})^{sur} = ({\mathcal A}^{min})^{sur} =
{\mathcal A}^{sur}\circ{\mathcal K}.$$
\item{(b$_2$)} If ${\mathcal A}$ is left accessible and also has the $l_1$-lifting
property, then we have
$$({\mathcal A}^{com})^{sur} = ({\mathcal A}^{min})^{sur} =
{\mathcal K}\circ{\mathcal A}^{sur}.$$}\par
\end{corollary}
\pf
Since $\mathcal A$ is left accessible, $\overline{\mathcal F}\circ{\mathcal A}
= {\mathcal A}^{min} \subset {\mathcal A}^{com}$. Then by Proposition 5.1(2$^\prime$),
we have
$${\mathcal K}\circ{\mathcal A}^{inj} = (\overline{\mathcal F}\circ{\mathcal A})^{inj}
 = ({\mathcal A}^{min})^{inj} \subset ({\mathcal A}^{com})^{inj}.$$
Since ${\mathcal A}^{com} \subset {\mathcal K}\circ{\mathcal A} \subset
{\mathcal K}\circ{\mathcal A}^{inj}$ and since ${\mathcal K}\circ{\mathcal A}^{inj}$
is injective, (a$_1$) follows. Dualising the arguments (b$_1$) also follows. The remaining
statements (a$_2$) and (b$_2$) follow directly from Propositions 5.1(c$^\prime$) and
5.2(c$^\prime$) respectively.
   \hfill $\triangle$
\vskip 8pt plus 0fill
The case of the injective-surjective hull is much simpler as we see now.
\begin{proposition}
$({\mathcal A}^{com})^{inj \hspace{1mm} sur} =  ({\mathcal A}^{min})^{inj \hspace{1mm} sur}$.
\end{proposition}
\pf
Let $T \in ({\mathcal A}^{com})^{inj \hspace{1mm} sur}(X, Y)$, for any Banach spaces
$X$ and $Y$. Then $i_Y\circ Tq_X \in {\mathcal A}^{com}(L_1(B_X), l_\infty (B_{Y^*}))$.
Thus there are Banach spaces
$Z$ and $W$ and operators $V \in {\mathcal K}(l_1(B_X), Z)$, $U \in {\mathcal K}
(W, l_\infty(B_{Y^*}))$ and $S \in {\mathcal A}(Z, W)$ such that $i_Y\circ T\circ q_X =
U\circ S\circ V$.
Since the Banach spaces $(l_1(B_X))^*$ and $l_\infty(B_{Y^*})$ have the approximation property,
we obtain $V \in {\overline{\mathcal F}}(l_1(B_X), Z)$ and $U \in {\overline{\mathcal F}}$.
Thus $i_Y\circ T\circ q_X \in ({\overline{\mathcal F}}\circ {\mathcal A}\circ{\overline{\mathcal F}})
(l_1(B_X), l_\infty(B_{Y^*}))$ and consequently, $T \in
({\mathcal A}^{min})^{inj \hspace{1mm} sur}$. Hence
$({\mathcal A}^{com})^{inj \hspace{1mm} sur} \subset
({\mathcal A}^{min})^{inj \hspace{1mm} sur}$. As ${\mathcal A}^{min}
\subset {\mathcal A}^{com}$, the result follows.
   \hfill $\triangle$
\vskip 8pt plus 0fill

\section{Approximating ${\mathcal A}\circ{\mathcal K}$ and ${\mathcal K}\circ{\mathcal A}$ by finite rank operators}
A Banach space $X$ is said to have the approximation property if given a
compact set $K\subset X$ and an $\epsilon >0$, there is a finite rank operator $T$ on $X$
such that $\|Tx-x\|<\epsilon$, for all $x\in K$. Recall that Grothendieck \cite{G} showed
that the following statements are equivalent:
\newline{\indent(1) The Banach space X has the approximation property.}
\newline{\indent(2) For every Banach space Y , the finite rank operators ${\mathcal F}(Y,X)$ is dense
in ${\mathcal B}(Y,X)$ in the topology of uniform convergence on compact sets.}
\newline{\indent(3) For every Banach space $Y$ , the finite rank operators ${\mathcal F}(X,Y)$ is dense
in ${\mathcal B}(X,Y)$ in the topology of uniform convergence on
compact sets.}
\newline{\indent(4) For every Banach space $Y$ , the finite rank operators ${\mathcal F}(Y,X)$ is dense
in compact operators ${\mathcal K}(Y,X)$ in the operator norm.}
\newline\noindent For the approximation property in the dual of a Banach space we have \cite[Theorem 1.e.5]{LT}
\newline{\indent(5) The dual $X^*$ of a Banach space X has the approximation property if
and only if for every Banach space $Y$ , the finite rank operators ${\mathcal F}(X,Y)$ is
dense in compact operators ${\mathcal K}(X,Y)$ in the operator norm.}

The authors in \cite[Definitions 2.1, 2.2 and 6.1 ]{SK1} introduced the notions of a $p$-compact set, a $p$-compact operator and
the corresponding notion of the $p$-approximation property, which is to approximate the identity operator by finite rank operators
on $p$-compact sets. However, the $p$-approximation property is equivalent to the density of finite rank operators in the
ideal of $p$-compact operators ${\mathcal K}_p$ in a norm weaker than the ideal norm of ${\mathcal K}_p$ \cite[Theorem 6.3]{SK1}.
Thus a prototype of the above equivalence ($1$)$\Leftrightarrow$($4$), for $1\leq p<\infty$, could not be achieved.  Later,
the authors \cite[Definition 4.4]{SK2}, taking a que from ($2$) above, introduced another notion, namely, the approximation property of type $p$, for $1\leq p \leq \infty$,
in terms of a locally convex topology $\lambda _p$ on the operator ideal $\Pi _p$. It was proved \cite[Theorems 4.5 and 4.6]{SK2} that both
($2$)$\Leftrightarrow$($4$) and ($5$) of the above list have suitable prototypes for the approximation property of type $p$. To be precise,
a Banach space $X$ has the approximation property of type $p$ if and only if the finite rank operators are dense in ${\mathcal K}_p^{dual}(Y,X)$ for
all Banach spaces $Y$, and the dual space $X^*$ has the approximation property of type $p$ if and only if the finite rank operators are dense
in ${\mathcal K}_p^{dual \ dual}(X,Y)$ for all Banach spaces $Y$.
In what follows, we seek to reinvent the later for a general operator ideal, and in particular, for injective and surjective operator ideals.

First, we record an observation due to Grothendieck, that is essentially contained in
his proof of the several equivalent formulations of the approximation property \cite[Theorem 1.e.4 $(4) \Rightarrow (1)$]{LT}.

\begin{lemma}
Let $(X, \|.\|)$ be a Banach space and $K$ be a compact subset of $X$. Then there is a Banach space
$(Y, \|.\|_0)$ be formally contained in $X$ such that the formal identity map $i_Y:Y\to X$ is a
compact operator with $K \subset \overline{i_Y(Ball Y)}$. Furthermore, given  a continuous function $g$
on $Y$ and $\epsilon >0$, there exists an $f\in X^*$ such that
$$|f(y) - g(y)|<\epsilon, for \ all \ y\in i_Y(Ball Y).$$
\end{lemma}

Next, again in the spirit of ($2$) above,  we define a locally convex topology on ${\mathcal A}(X, Y)$ for an operator ideal $\mathcal A$
followed by a corresponding approximation property in the following manner.

\begin{definition}
Let $\mathcal A$ be a operator ideal.
For Banach spaces $X$, $Y$ and a compact set $K \subset X$ we define a seminorm $\|.\|_K$ on
${\mathcal A}(X, Y)$ given by
$$\|T\|_K=\inf\{\alpha _{\kappa}(T\circ i_Z): i_Z: Z\to X \ as \ above \},$$
where $\alpha _{\kappa}(T) =\inf\{ \alpha(T_1) \|T_2\|: T=T_1\circ T_2\in{\mathcal A}\circ{\mathcal K}\}$
is the usual composition norm on ${\mathcal A}\circ{\mathcal K}$. Then the family of seminorms
$\{\|.\|_K: K\subset X \ is \ compact\}$ defines a locally convex topology $\lambda _{\mathcal A}$ on
${\mathcal A}\circ{\mathcal K}(X, Y)$.
\end{definition}

\begin{definition}
Let $\mathcal A$ be a operator ideal. A Banach spaces $X$ is said to have the $\mathcal A$-approximation
property ($\mathcal A$-a.p., for short) if for every Banach space $Y$, ${\mathcal F}(Y, X)$ is dense in
${\mathcal A}(Y, X)$ in the $\lambda _{\mathcal A}$-topology.
\end{definition}

Note that for the operator ideal $\mathcal B$, the topology $\lambda _{\mathcal B}$ is the topology on
${\mathcal B}(X, Y)$ of uniform convergence on compact sets of $X$. Recently, the authors \cite[Definition 4.3]{SK2}
have defined for each $1\leq p \leq \infty$, a locally convex topology $\lambda _p$ on $\Pi _p(X, Y)$. In terms of
the Definition 6.2 above  this is the $\Pi _p$-topology on $\Pi _p(X, Y)$. A Banach space $X$ is said to have the
approximation property of {\it type} $p$ \cite[Definition 4.4]{SK2} if for any Banach space $Y$, the finite
rank operators ${\mathcal F}(Y, X)$ is dense in $\Pi _p(Y, X)$ in the $\lambda _p$-topology. Thus in terms of
the Definition 6.3 above this is the $\Pi _p$-approximation property on $X$.

Now, we characterise the $\mathcal A$-approximation property in terms of the density of finite rank operators
$\mathcal F$ in ${\mathcal A}\circ{\mathcal K}$ in its composition operator norm $\alpha _{\kappa}$.

\begin{theorem}
Let $\mathcal A$ be an operator ideal. A Banach space $X$ has the $\mathcal A$-approximation property if and only
if for every Banach space $Y$, we have
$${\mathcal A}\circ{\mathcal K}(Y, X)\subset\overline{{\mathcal F}(Y, X)}^{\alpha _{\kappa}}.$$
\end{theorem}

\pf
Let $T\in {\mathcal A}\circ{\mathcal K}(Y, X)$ and $\epsilon >0$. Consider a factorisation $T=T_1\circ T_2$ where
$T_1\in{\mathcal A}(G,X)$ and $T_2\in{\mathcal K}(Y,G)$ for some Banach space $G$. Set $K=\overline{T_2(Ball Y)}$.
Since $X$ has the $\mathcal A$-a.p., there is a $S_0\in{\mathcal F}(G,X)$ such that
$$\|T_1-S_0\|<\epsilon /\|T_2\|.$$
Furthermore, as in the above Lemma 6.1, we can find a Banach space $Z$ formally contained in $G$ such that $i_Z:Z\to G$
is compact and we have
$$\alpha _{\kappa}((T_1-S_0)\circ i_Z)< \epsilon /\|T_2\|.$$
Since $i_Z(k)=k$ for all $k\in K=\overline{T_2(Ball Y)}$, we get $i_Z\circ T_2=T_2$. Set $S=S_0\circ T_2$. Then
$S\in {\mathcal F}(Y,X)$ and we have
$$\alpha _{\kappa}(T-S) = \alpha _{\kappa}(T_1\circ T_2-S_0\circ T_2)\leq \alpha _{\kappa}((T_1-S_0)i_Z)\|T_2\|<\epsilon.$$

Conversely, let $T\in{\mathcal A}(Y,X)$, $K\subset Y$ a compact set and $\epsilon >0$. By Lemma 6.1, there is a Banach
space $Z$ formally contained in $Y$ such that $i_Z:Z\to Y$ is compact and $K\subset i_Z(Ball Z)$. Then
$T\circ i_Z\in{\mathcal A}\circ{\mathcal K}(Z,X)$ so that there is an $S_1\in {\mathcal F}(Z,X)$ such that
$\alpha _{\kappa}(T\circ i_Z-S_1)<\epsilon /2$. Consider a factorisation $S_1=S_2\circ S_3$ for some operators
$S_2\in{\mathcal F}(G,X)$, $S_3\in{\mathcal F}(Z,G)$ and some Banach space $G$. Choose $S_4\in{\mathcal F}(Z,G)$ such
that $\|S_3-S_4\|<\epsilon /4\alpha _{\kappa}(S_2)$. Let $S_4 \sim\sum _{i=1}^{n}{f_i\otimes x_i}$, with
$f_1, f_2,\cdots f_n\in Z^*$ and $x_1, x_2,\cdots x_n\in G$. Again by Lemma 6.1, given $\delta >0$ we can find
$g_1, g_2,\cdots g_n\in Y^*$ such that $|f_i(z)-g_i(z)|<\delta$ for all $z\in Ball Z$. Indeed, we may choose
$\delta =\epsilon /4n\alpha _{\kappa}(S_2)\max_{1\leq i\leq n}\|x_i\|$ such that $\|S_4-S^\prime\|<\epsilon /4\alpha _{\kappa}(S_2)$,
where $S^\prime=\sum _{i=1}^n{g_i\cdot i_Z\otimes x_i}$. Set $S^{\prime\prime}=\sum _{i=1}^n{g_i\otimes x_i}$ so that
$S^{\prime\prime}\in {\mathcal F}(Y,G)$ and $S^{\prime}= S^{\prime\prime}\cdot i_Z$.
Thus $S=S_2\cdot S^{\prime\prime}\in {\mathcal F}(Y,X)$ and
$$\alpha _{\kappa}((T-S)\cdot i_Z\leq\alpha _{\kappa}(T\cdot i_Z-S_1)+\alpha _{\kappa}(S_1-S_2\cdot S^{\prime\prime}\cdot i_Z)$$
$$\leq\epsilon /2+\alpha _{\kappa}(S_2\cdot S_3-S_2\cdot S^{\prime})$$
$$\leq\epsilon /2+\alpha _{\kappa}(S_2)\|S_3-S^{\prime}\|$$
$$\leq\epsilon /2+\alpha _{\kappa}(S_2)\{\|S_3-S_4\|+\|S_4-S^{\prime}\|\}<\epsilon .$$
It follows that $\|T-S\|_K\leq\alpha _{\kappa}((T-S)\cdot i_Z)<\epsilon$, so that $X$ has the $\mathcal A$-a.p.
 \hfill $\triangle$

\vskip 8pt plus 0fill
It was proved by the authors \cite[Theorem 4.5]{SK2} that a Banach space $X$ has the approximation property of {\it type} $p$
if and only if for every Banach space $Y$, the finite rank operators ${\mathcal F}(Y,X)$ is dense in the composition operator
ideal $\Pi _p\cdot{\mathcal K}(Y,X)$ in its factorisation norm ${\pi _p}_{\kappa}$. Thus, Theorem 4.5 in \cite{SK2} is a special case
of Theorem 6.4 above, for the operator ideal $\Pi _p$.

\begin{theorem}
Let $\mathcal A$ be an operator ideal and $X$ a Banach space. If the dual space $X^*$ has the $\mathcal A$-
approximation property then for every Banach space $Y$ we have
$${\mathcal K}\cdot{\mathcal A}^{dual}(X,Y)\subset\overline{{\mathcal F}(X,Y)}^{{_\kappa}{\alpha ^d}}.$$
The converse holds if, in addition, ${\mathcal A}\subset{\mathcal A}^{dual \ dual}$ (equivalently,
${\mathcal A}^{dual \ dual}={\mathcal A}^{reg})$.
\end{theorem}
\pf
Let $X^*$ have the $\mathcal A$-a.p. Then by Theorem 6.4 above, for every Banach space $Y$, we have
$${\mathcal A}\circ{\mathcal K}(Y,X^*)\subset\overline{{\mathcal F}(Y,X^*)}^{\alpha _{\kappa}}.$$
It follows that ${\mathcal A}\circ{\mathcal K}(Y^*,X^*)\subset\overline{{\mathcal F}(Y^*,X^*)}^{\alpha _{\kappa}}$, for all
Banach spaces $Y$. Let $T\in{\mathcal K}\cdot{\mathcal A}^{dual}(X,Y)$ and $\epsilon >0$. Then there is a Banach
space $G$ and operators $T_1\in{\mathcal K}(G,Y)$ and $T_2\in{\mathcal A}^{dual}X,G)$ such that $T=T_1\circ T_2$.
Then $T^*=T_2^*\circ T_1^*\in{\mathcal A}\circ{\mathcal K}(Y^*,X^*)$ so that there is an $S\in{\mathcal F}(Y^*,X^*)$
with $\alpha _{\kappa}(T^*-S)<\epsilon$. We may write $T^*-S=S_1\circ S_2$, where $S_1\in{\mathcal A}(Z,X^*)$ and
$S_2\in{\mathcal K}(Y^*,Z)$ for some Banach space $Z$ such that $\alpha(S_1)<\epsilon$ and $\|S_1\|=1$. Now, $S^*\in{\mathcal F}(X^{**},Y)$ and
$${_\kappa}{\alpha ^d}(T-S^*\circ \sigma _X)= {_\kappa}\alpha ^d(T^{**}\circ\sigma _X-S^*\circ\sigma _X)$$
$$\leq {_\kappa}\alpha ^d(T^{**}-S^*)\|\sigma _X\|)$$
$$\leq\|S_2^*\|\alpha ^d(S_1^*)$$
$$=\alpha   ^{**}(S_1)\leq\alpha(S_1)<\epsilon ,$$
where $\sigma _X:X\to X^{**}$ is the canonical embedding.

Conversely, let $\mathcal A$ be an operator ideal satisfying  ${\mathcal A}\subset{\mathcal A}^{dual \ dual}$,
$T\in{\mathcal A}\circ{\mathcal K}(Y,X^*)$ and $\epsilon >0$. Then there is a Banach space $G$ and operators $T_1\in{\mathcal A}(G,X^*)$ and
$T_2\in{\mathcal K}(Y,G)$ such that $T=T_1\circ T_2$. It follows that $T^*\in {\mathcal K}\circ
{\mathcal A}^{dual}(X^{**},Y^*)$, so that $T^*\circ\sigma _X\in{\mathcal K}\circ{\mathcal A}^{dual}(X,Y^*)$, and by our assumption there
is an $S\in\overline{\mathcal F}(X,Y^*)$ such that ${_\kappa}{\alpha ^d}(T^*\circ\sigma _X-S)<\epsilon$. Now we can choose a Banach space
$Z$ and operators $S_1\in{\mathcal K}(Z,Y^*)$ an $S_2\in{\mathcal A}^{dual}(X,Z)$ such that $T^*\circ\sigma_X-S=S_1\circ S_2$ with $\|S_1\|=1$
and $\alpha ^d(S_2)<\epsilon$. Thus $S^*\in\overline{\mathcal F}(Y^{**},X^*)$ and we have
$$\alpha _{\kappa}(T-S^*\circ\sigma _Y)=\alpha _{\kappa}(\sigma _X^*T^{**}\circ\sigma _Y-S^*\sigma _Y)$$
$$=\alpha _{\kappa}(S_2^*\circ S_1^*\circ\sigma _Y)$$
$$\leq\alpha(S_1^*)\|S_1^*\|\|\sigma _Y\|$$
$$=\alpha ^d(S_2)<\epsilon.$$
It follows that $X^*$ has the $\mathcal A$-a.p.   \hfill $\triangle$

The authors in \cite[Theorem 4.6]{SK2} have shown that the dual $X^*$ of a Banach space $X$ has the approximation property of {\it type} $p$
if and only if ${\mathcal F}(X,Y)$ is dense in ${\mathcal K}\cdot{\Pi _p^{dual}}(X,Y)$ in its composition ideal norm for all Banach spaces$Y$.
Thus, Theorem 4.6 in \cite{SK2} is a special case of Theorem 6.5 above.

Finally, at the end we discuss certain special cases of the $\mathcal A$-a.p. for injective and surjective operator ideals $\mathcal A$ with
$\ell _{\infty}$-extension or $\ell _1$-lifting property.

\begin{corollary}
Let $\mathcal A$ be a injective operator ideal with the $\ell _{\infty}$-extension property and $X$ a Banach space. Then we have
the following: \par
{\leftskip=6mm
\item{(a)} $X$ has the $\mathcal A$-approximation property if and only if
$$({\mathcal A}^{min})^{inj}(Y,X)={\mathcal A}^{min}(Y,X),$$
for all Banach spaces $Y$.
\item{(b)} The dual space $X^*$ has the $\mathcal A$-approximation property if and only if
$$(({\mathcal A}^{dual})^{min})^{sur}(X,Y)=({\mathcal A}^{dual})^{min}(X,Y),$$
for all Banach spaces $Y$. \par}
\end{corollary}

\pf
(a) Since $\mathcal A$ is injective, it is right accessible and it also has the $\ell _\infty$-extension property. Thus, by Corollary 5.3($a_2$)
we have $({\mathcal A}^{min})^{inj}={\mathcal A}^{inj}\circ{\mathcal K}={\mathcal A}\circ{\mathcal K}$.
By Theorem 6.4, $X$ has the $\mathcal A$-a.p. if and only if ${\mathcal A}^{min}(Y,X)$ is dense in ${\mathcal A}\circ{\mathcal K}(Y,X)$
for all Banach spaces $Y$.  Thus $X$ has the $\mathcal A$-a.p. if and only if $({\mathcal A}^{min})^{inj}(Y,X)={\mathcal A}^{min}(Y,X)$, for all Banach spaces $Y$.

(b) Since $\mathcal A^{dual}$ is surjective, it is left asseccible and by Theorem 4.9 it also has the $\ell _1$-lifting property. Thus by Corollary 5.3($b_2$) we have
$(({\mathcal A}^{dual})^{min})^{sur}=({\mathcal A}^{dual})^{min}$. Now, by Theorem 6.5 , $X^*$ has the $\mathcal A$-a.p.
if and only if $({\mathcal A}^{dual})^{min}(X,Y)$ is dense in ${\mathcal K}\circ{\mathcal A}^{dual}(X,Y)$ for all Banach spaces $Y$. Thus $X^*$ has the
$\mathcal A$-a.p. if and only if $(({\mathcal A}^{dual})^{min})^{sur}(X,Y)=({\mathcal A}^{dual})^{min}(X,Y)$, for all Banach spaces $Y$.  \hfill $\triangle$

Again note that since $\Pi _p$ is injective and has the $\ell _\infty$-extension property, Theorems 4.5 and 4.6 in \cite{SK2} are special cases of the above corollary.
Let us also note here that we do not know whether there is an injective operator ideal that fails the $\ell _\infty$-extension property.
\begin{corollary}
Let $\mathcal A$ be a right accessible surjective operator ideal and $X$ be a Banach space. Then we have the following: \par
{\leftskip=6mm
\item{(a)} $X$ has the $\mathcal A$-approximation property if and only if
$$({\mathcal A}^{min})^{sur}(Y,X)={\mathcal A}^{min}(Y,X),$$
for all Banach spaces $Y$.
\item{(b)} The dual space
$X^*$ has the $\mathcal A$-approximation property if and only if
$$(({\mathcal A}^{dual})^{min})^{inj}(X,Y)=({\mathcal A}^{dual})^{min}(X,Y),$$
for all Banach spaces $Y$. \par}
\end{corollary}

\pf
(a) By Corollary 5.3($b_1$), we have $({\mathcal A}^{min})^{sur}={\mathcal A}\circ{\mathcal K}$. Now it follows by Theorem 6.4 that $X$ has the $\mathcal A$-a.p.
if and only if $({\mathcal A}^{min})^{sur}(Y,X)={\mathcal A}^{min}(Y,X)$, for all Banach spaces $Y$.

(b) Since ${\mathcal A}^{dual}$ is injective, by Corollary 5.3($a_1$), we have $(({\mathcal A}^{dual})^{min})^{inj}={\mathcal K}\circ{\mathcal A}^{dual}$. Now it follows by
Theorem 6.5 that $X^*$ has the $\mathcal A$-a.p. if and only if $(({\mathcal A}^{dual})^{min})^{inj}(X,Y)=({\mathcal A}^{dual})^{min}(X,Y)$, for all Banach spaces $Y$.
 \hfill $\triangle$
\vskip 8pt plus 0fill


\begin{thebibliography}{150}

\bibitem{DF} A.\ Defant and K.\ Floret, {\it Tensor norms and operator ideals},
North-Holland, Amsterdam, 1993.

\bibitem{DJT} J.\ Diestel , H.\ Jarchow and A.\ Tonge, {\it Absolutely Summing Operators},
Cambridge Univ.\  Press,  1995.

\bibitem{G} A.\ Grothendieck, {\it Produits tensoriels topologiques et especas nucleaires},
Memo.\ Amer.\ Math.\ Soc.\ {\bf 16}(1955).

\bibitem{LT} J.\ Lindenstrauss and L.\ Tzafriri, {Classical Banach spaces - Sequence spaces},
Springer-Verlag, Berlin, 1977.

\bibitem{Or} F.\ Ortel, {\it Composition of operator ideals and their regular hulls},
Acta.\ Univ.\ Carolin.\ Math.\ Phys.\ {\bf 36}(1995), 69-72.

\bibitem{P1} A.\ Pietsch, {\it Operator Ideals}, North Holland, Amsterdam, 1980.

\bibitem{SK1} D.\ P.\ Sinha and A.\ K.\ Karn, {\it Compact operators whose adjoints factor
through subspaces of $l_p$}, Studia Math.\, {\bf 150}(2002), 17-33.

\bibitem{SK2} D.\ P.\ Sinha and A.\ K.\ Karn, {\it Compact operators which factor
through subspaces of $l_p$}, Math.\ Nachr.\, {\bf 281}(2008), 412-423.

\end{thebibliography}
\end{document}